\documentclass[12pt]{article}
\usepackage{latexsym}
 \usepackage{mathrsfs}
\usepackage{amssymb}
\usepackage{graphicx}
\usepackage{cite}

\newtheorem{Theorem}{Theorem}[part]
\newtheorem{Definition}{Definition}[part]
\newtheorem{Proposition}{Proposition}[part]

\newtheorem{Lemma}{Lemma}[part]

\newtheorem{Example}{Example}[part]

\def \ep{\hbox{ }\hfill$\Box$}

\begin{document}
\title{On $Q$-Tensors \thanks{This work was partially supported by the National Natural
Science Foundation of China (Grant No. 11431002).}
}

\author{
Zheng-Hai Huang\thanks{Corresponding Author. Department of Mathematics, School of Science, Tianjin University, Tianjin 300072, P.R. China. This author is also with the Center for Applied Mathematics of Tianjin University. Email: huangzhenghai@tju.edu.cn. Tel:+86-22-27403615 Fax:+86-22-27403615}
\and
Yun-Yang Suo\thanks{Department of Mathematics, School of Science, Tianjin University, Tianjin 300072, P.R. China. Email: suoyunyang@163.com}
\and
Jie Wang\thanks{Department of Mathematics, School of Science, Tianjin University, Tianjin 300072, P.R. China. Email: wxinnuan@163.com}
}

\date{September 10, 2015}

\maketitle

\begin{abstract}
\noindent

One of the central problems in the theory of linear complementarity problems (LCPs) is to study the class of $Q$-matrices since it characterizes the solvability of LCP. Recently, the concept of $Q$-matrix has been extended to the case of tensor, called $Q$-tensor, which characterizes the solvability of the corresponding tensor complementarity problem -- a generalization of LCP; and some basic results related to $Q$-tensors have been obtained in the literature. In this paper, we extend two famous results related to $Q$-matrices to the tensor space, i.e., we show that within the class of strong $P_0$-tensors or nonnegative tensors, four classes of tensors, i.e., $R_0$-tensors, $R$-tensors, $ER$-tensors and $Q$-tensors, are all equivalent. We also construct several examples to show that three famous results related to $Q$-matrices cannot be extended to the tensor space; and one of which gives a negative answer to a question raised recently by Song and Qi.
\vspace{3mm}

\noindent {\bf Key words:}\hspace{2mm} $Q$-tensor, tensor complementarity problem, strong $P_0$-tensor, nonnegative tensor \vspace{3mm}

\noindent {\bf Mathematics Subject Classifications (2000):}\hspace{2mm} 90C33, 65K10, 15A18, 15A69, 65F15, 65F10. \vspace{3mm}

\end{abstract}
\section{Introduction}
\setcounter{equation}{0} \setcounter{Assumption}{0}
\setcounter{Theorem}{0} \setcounter{Proposition}{0}
\setcounter{Corollary}{0} \setcounter{Lemma}{0}
\setcounter{Definition}{0} \setcounter{Remark}{0}
\setcounter{Algorithm}{0}

\hspace{4mm} For any given $A\in \mathbb{R}^{n\times n}$ and $q\in \mathbb{R}^n$, the linear complementarity problem \cite{cps-92}, denoted by LCP$(q,A)$, is to find a vector $x\in \mathbb{R}^n$ such that
\begin{eqnarray*}
x\geq 0,\quad Ax+q\geq 0,\quad x^T(Ax+q)=0.
\end{eqnarray*}
If LCP$(q,A)$ has a solution for every vector $q\in \mathbb{R}^n$, we say that $A$ is a {\it $Q$-matrix} \cite{murty-72}. Since the solvability of LCP$(q,A)$ is very important, much research in the theory of LCP$(q,A)$ has been devoted to finding constructive characterizations of subclasses of $Q$-matrices; and fruitful results have been obtained \cite{ac-79,Danao-94,Gowda-90,jp-85,jp-89,mpr-95,Pang-79}.

A real $m$-order $n$-dimensional tensor can be denoted by $\mathscr{A}=(a_{i_1\cdots i_m})$ with $a_{i_1\cdots i_m}\in \mathbb{R}$ for all $i_j\in \{1,2,\ldots,n\}$ and $j\in \{1,2,\ldots,m\}$. Obviously, $\mathscr{A}$ is a matrix when $m=2$.
As an extension of matrix, tensor has been widely studied (see an excellent
survey by Kolda and Bader \cite{kb-09} and references therein).
For any given real $m$-order $n$-dimensional tensor $\mathscr{A}$ and vector $q\in \mathbb{R}^n$, the tensor complementarity problem \cite{bhw-15-r,dlq-2015-r,sq-15-r,sq-15,sq-15-r1,whb-15-r}, denoted by TCP$(q,\mathscr{A})$, is to find a vector $x\in \mathbb{R}^n$ such that
\begin{eqnarray*}
x\geq 0,\quad \mathscr{A}x^{m-1}+q\geq 0,\quad x^T(\mathscr{A}x^{m-1}+q)=0,
\end{eqnarray*}
where $\mathscr{A}x^{m-1}\in \mathbb{R}^n$ with
$$
\left(\mathscr{A}x^{m-1}\right)_i:=\sum_{i_2,\ldots,i_m=1}^na_{ii_2\cdots i_m}x_{i_2}\cdots x_{i_m},\;\forall i\in \{1,2,\ldots,n\}.
$$
Recently, Song and Qi \cite{sq-15-r} extended the concept of $Q$-matrix to the case of tensor, called $Q$-tensor, i.e., a real $m$-order $n$-dimensional tensor $\mathscr{A}$ is called a {\it $Q$-tensor} if TCP$(q,\mathscr{A})$ has a solution for every vector $q\in \mathbb{R}^n$. Furthermore, they proved that several classes of tensors are the subclasses of $Q$-tensors.

In this paper, we consider several famous results related to $Q$-matrices and investigate whether they can be extended to the tensor space or not. In Section 2, after reviewing some basic concepts and known related results, we give a sufficient condition to judge whether a tensor is a $Q$-tensor or not.

Recall that a matrix $A\in \mathbb{R}^{n\times n}$ is called a {\it $P_0$-matrix} if all its principal minors are nonnegative; an {\it $R_0$-matrix} if LCP$(0,A)$ has a unique solution; and an {\it $R$-matrix} if LCP$(te,A)$ has a unique solution for each scalar $t\geq 0$, where $e\in \mathbb{R}^n$ denotes the vector of all ones. In 1979, Agangic and Cottle \cite{ac-79} proved that within the class of $P_0$-matrices, both the classes of $R$-matrices and $R_0$-matrices are equivalent to the class of $Q$-matrices. Since the classes of matrices mentioned above play important roles in the field of variational inequalities and complementarity problems \cite{fp-03,hp-90,hxq-06,hqs-04,kanzow-96,kmn-91}, it is important to extend Agangic-Cottle's result to the tensor space. It is known that $P_0$-matrix, $R_0$-matrix and $R$-matrix have been extended to the case of tensor by Song and Qi \cite{sq-15,sq-15-r}, called $P_0$-tensor, $R_0$-tensor and $R$-tensor, respectively. Recently, Ding, Luo and Qi \cite{dlq-2015-r} gave another extension of $P_0$-matrix, which is called $P_0^\prime$-tensor in this paper. A natural question is whether Agangic-Cottle's result can be extended to the tensor space or not. In Section 3, we answer this question. We clarify that within the class of $P_0$-tensors (or $P_0^\prime$-tensors), the above equivalence cannot be extended to the tensor space. In order to extend Agangic-Cottle's result to the tensor space, we introduce a new class of tensors, called strong $P_0$-tensors, which is a generalization of $P_0$-matrices; and show that within the class of strong $P_0$-tensors, the classes of $R_0$-tensors, $R$-tensors and $ER$-tensors are all equivalent to the class of $Q$-tensors, where the class of $ER$-tensors was recently introduced by Wang, Huang and Bai \cite{whb-15-r}. In this section, we also discuss the relationships among $P_0$-tensors, $P_0^\prime$-tensors and strong $P_0$-tensors.

Recall that a matrix $A\in \mathbb{R}^{n\times n}$ is said to be {\it nonnegative} if all its elements are nonnegative; and a real $m$-order $n$-dimensional tensor $\mathscr{A}$ is said to be {\it nonnegative} if all its elements are nonnegative. In 1994, Danao \cite{Danao-94} proved that within the class of nonnegative matrices, the class of $Q$-matrices coincides with the class of $R_0$-matrices.
Since many tensors from practical problems are nonnegative, nonnegative tensors have been extensively studied in the last ten years \cite{cpz-08,cpz-11,fgh-13,hhq-14,nqz-09,yy-10}. It is interesting whether Danao's result can be extended to the tensor space or not. In Section 4, we answer this question. We show that within the class of nonnegative tensors, the classes of $R_0$-tensors, $R$-tensors and $ER$-tensors are all equivalent to the class of $Q$-tensors. In addition, we also discuss the relationship between the class of nonnegative tensors and the  class of strong $P_0$-tensors.

Recall that a matrix $A\in \mathbb{R}^{n\times n}$ is said to be {\it semi-monotone} if for every $x\in \mathbb{R}^n$ with $0\neq x\geq 0$, there exists an index $i$ such that $x_i>0$ and $(Ax)_i\geq 0$; and {\it copositive} if $x^TAx\geq 0$ for all $x\in \mathbb{R}^n$ satisfying $x\geq 0$. In Section 5, we consider two results obtained by Pang \cite{Pang-79} which are related to semi-monotone matrices and $Q$-matrices; and a result obtained by Jeter and Pye \cite{jp-85} which is related to copositive matrices and $Q$-matrices. We illustrate that these three results cannot be extended to the tensor space by using several examples. In addition, the final conclusions are given in Section 6.

In the rest of this paper, we assume that $m\geq 3$ and $n\geq 2$ are two integers unless otherwise specialized; and use $\mathbb{T}_{m,n}$ to denote the set of all real $m$-order $n$-dimensional tensors.
For any positive integer $n$, we denote $[n]:=\{1,2,\ldots,n\}$ and $\mathbb{R}^n_+:=\{x\in \mathbb{R}^n: x\geq 0\}$. We use $\mathbb{T}_{m,n}\bigcap Q$ to denote the set of all real $m$-order $n$-dimensional $Q$-tensors and $\mathbb{R}^{n\times n}\bigcap Q$ to denote the set of all real $Q$-matrices of $n\times n$. A tensor $\mathscr{A}\in \mathbb{T}_{m,n}$ satisfying $\mathscr{A}\in Q$ (i.e., $\mathscr{A}\in \mathbb{T}_{m,n}\bigcap Q$) means that $\mathscr{A}$ is a real $m$-order $n$-dimensional $Q$-tensor and a matrix $A\in \mathbb{R}^{n\times n}$ satisfying $A\in Q$ (i.e., $A\in \mathbb{R}^{n\times n}\bigcap Q$) means that $A$ is a real $Q$-matrix of $n\times n$. Similar notation is used for other classes of tensors or matrices.

\section{Preliminaries}
\setcounter{equation}{0} \setcounter{Assumption}{0}
\setcounter{Theorem}{0} \setcounter{Proposition}{0}
\setcounter{Corollary}{0} \setcounter{Lemma}{0}
\setcounter{Definition}{0} \setcounter{Remark}{0}
\setcounter{Algorithm}{0}

\hspace{4mm}
In this section, we review definitions and properties of several structured tensors, which are useful for our subsequent discussions. We also give a sufficient condition to judge whether a tensor is a $Q$-tensor or not.

\begin{Definition} \label{p0def}
A tensor $\mathscr{A}\in \mathbb{T}_{m,n}$ is said to be
\begin{itemize}
  \item[(i)] {\bf a $P_0$-tensor} iff for each $x\in \mathbb{R}^n\setminus\{0\}$, there exists an index $i\in [n]$ such that
      $$ x_i\neq 0\quad \mbox{\rm and}\quad x_i(\mathscr{A}x^{m-1})_i\geq0;$$
  \item[(ii)] {\bf a $P_0^\prime$-tensor} iff for each $x\in \mathbb{R}^n\setminus\{0\}$, there exists an index $i\in [n]$ such that
      $$x_i\neq 0\quad \mbox{\rm and}\quad x_i^{m-1}(\mathscr{A}x^{m-1})_i\geq0.$$
\end{itemize}
\end{Definition}

The concept of $P_0$-tensor was introduced by Song and Qi in \cite{sq-15}; and the concept of  $P_0^\prime$-tensor was introduced by Ding, Luo and Qi in \cite{dlq-2015-r} with the name of $P_0$-tensor. Since it is different from  $P_0$-tensor defined by Song and Qi in \cite{sq-15} (see the next section), in this paper, we call $P_0$-tensor introduced by Ding, Luo and Qi \cite{dlq-2015-r} as $P_0^\prime$-tensor to avoid confusion.
When $m=2$, a $P_0$-tensor or a $P_0^\prime$-tensor reduces to a $P_0$-matrix.

\begin{Definition}\label{rdef}
(i) A tensor $\mathscr{A}\in \mathbb{T}_{m,n}$ is called an $R$-tensor, if there exists no $(x, t)\in(\mathbb{R}^n_+\setminus \{0\})\times\mathbb{R}_+$ such that
\begin{equation}\label{rtensor}
\left\{\begin{array}{ll}
(\mathscr{A}x^{m-1})_i+t=0,\quad & \mbox{\rm if}\; x_i>0,\\
(\mathscr{A}x^{m-1})_i+t\geq0,\quad & \mbox{\rm if}\; x_i=0.
\end{array}\right.
\end{equation}
(ii) A tensor $\mathscr{A}\in \mathbb{T}_{m,n}$ is called an $R_0$-tensor, if the system (\ref{rtensor}) has no solution when $t=0$, i.e., there exists no $x\in \mathbb{R}^n_+\setminus \{0\}$ such that
$$
\left\{\begin{array}{ll}
(\mathscr{A}x^{m-1})_i=0,\quad & \mbox{\rm if}\; x_i>0,\\
(\mathscr{A}x^{m-1})_i\geq0,\quad & \mbox{\rm if}\; x_i=0.
\end{array}\right.
$$
(iii) A tensor $\mathscr{A}\in \mathbb{T}_{m,n}$ is called an ER-tensor, if there exists no $(x, t)\in(\mathbb{R}^n_+\setminus \{0\})\times\mathbb{R}_+$ such that
$$
\left\{\begin{array}{lcll}
(\mathscr{A}x^{m-1})_i+tx_i&=&0,   \quad & \mbox{\rm if}\; x_i>0,\\
(\mathscr{A}x^{m-1})_i     &\geq&0,\quad & \mbox{\rm if}\; x_i=0.
\end{array}\right.
$$
\end{Definition}

The concepts of $R$-tensor and $R_0$-tensor were introduced by Song and Qi in \cite{sq-15-r}; and the concept of $ER$-tensor was introduced by Wang, Huang and Bai in \cite{whb-15-r}. It is obvious that an $R$-tensor or $ER$-tensor is an $R_0$-tensor.
When $m=2$, an $R$-tensor ($R_0$-tensor) reduces to an $R$-matrix ($R_0$-matrix) \cite{cps-92,K-72}; and an $ER$-tensor reduces to an $ER$-matrix \cite{whb-15-r}.

\begin{Definition} \label{semidef}\cite{sq-15-r}
A tensor $\mathscr{A}\in \mathbb{T}_{m,n}$ is said to be {\bf semi-positive} iff for each $x\in \mathbb{R}^n_+\setminus \{0\}$, there exists an index $i\in [n]$ such that $x_i>0$ and $(\mathscr{A}x^{m-1})_i\geq 0$.
\end{Definition}

Clearly, every $P_0$-tensor ($P_0^\prime$-tensor) is certainly semi-positive. It is shown that the class of semi-positive $R_0$-tensors is a subclass of $Q$-tensors \cite{sq-15-r}.
When $m=2$, a semi-positive tensor reduces to a semi-monotone matrix; and the set of all semi-monotone matrices is denoted by $L_1$ (or $E_0$) \cite{E-71,fp-66}.

\begin{Definition} \label{copdef}
A tensor $\mathscr{A}\in \mathbb{T}_{m,n}$ is said to be {\bf copositive} iff $\mathscr{A}x^m\geq 0$ for all $x\in \mathbb{R}^n_+$.
\end{Definition}

When $\mathscr{A}\in \mathbb{T}_{m,n}$ is symmetric, such a concept was first introduced by Qi \cite{Q-13}. When $m=2$, every copositive tensor reduces to a copositive matrix \cite{cps-92}.

\begin{Proposition}\label{pro1}
Suppose that $\mathscr{A}\in \mathbb{T}_{m,n}$. Then the following results hold.
\begin{itemize}
  \item [(i)]  If $\mathscr{A}\in R$, then $\mathscr{A}\in Q$.
  \item [(ii)] If $\mathscr{A}$ is semi-positive, then
  $$
  \mathscr{A}\in R_0 \quad \Longleftrightarrow \quad \mathscr{A}\in ER \quad \Longleftrightarrow \quad \mathscr{A}\in R.
  $$
  \item [(iii)] If $\mathscr{A}\in ER$, then $\mathscr{A}\in Q$.
  \item[(iv)] If $\mathscr{A}$ is nonegative, then $\mathscr{A}\in Q$ iff $a_{i\cdots i}>0$ for all $i\in [n]$.
\end{itemize}
\end{Proposition}

In Proposition \ref{pro1}, the results (i)-(iv) come from \cite[Theorem 3.2]{sq-15-r}, \cite[Theorem 3.3]{whb-15-r}, \cite[Corollary 4.1]{whb-15-r}, and \cite[Theorem 3.5]{sq-15-r}, respectively.

At the end of this section, we give a characterization of $Q$-tensor, which is an extension of Proposition 2.1 in \cite{mpr-95}.

\begin{Theorem}\label{thm-q-1}
Let $\mathscr{A}=(a_{i_{1}\cdots i_{m}})\in \mathbb{T}_{m,n}$. Denote $\mathscr{A}_{1\cdots1}:=(a_{i_{1}\cdots i_{m}})$ with $i_{1},\ldots,i_{m}\in [n]\setminus\{1\}$ and $\mathscr{A}_{2\cdots2}:=(a_{i_{1}\cdots i_{m}})$ with $i_{1},\ldots,i_{m}\in [n]\setminus\{2\}$. Suppose that $a_{1i_{2}\cdots i_{m}}=a_{2i_{2}\cdots i_{m}}$ for all $i_{2},\ldots, i_{m}\in [n]$, and both $\mathscr{A}_{1\cdots1}$ and $\mathscr{A}_{2\cdots2}$ are $Q$-tensors. Then $\mathscr{A}$ is a $Q$-tensor.
\end{Theorem}
{\bf Proof.} For any $q=(q_{1},\ldots ,q_{n})^{T}\in \mathbb{R}^n$, we consider the following two cases.

Case 1. $q_{2}\leq q_{1}$. In this case, we denote
$$
{\cal N}:=\left\{(i_{2},\ldots ,i_{m}): i_{2},\ldots ,i_{m}\in [n]\setminus\{1\}\right\}\quad \mbox{\rm and}\quad q_{-1}:=(q_{2},\ldots ,q_{n})^{T}.
$$
Then, for any $x=(0,\hat{x}^T)^T\in \mathbb{R}\times \mathbb{R}^{n-1}$ with $\hat{x}:=(x_2,\ldots,x_n)^T\in \mathbb{R}^{n-1}$, it follows that $\mathscr{A}_{1\cdots1}\hat{x}^{m-1}=((\mathscr{A}_{1\cdots1}\hat{x}^{m-1})_1,\ldots,(\mathscr{A}_{1\cdots1}\hat{x}^{m-1})_{n-1})^T\in \mathbb{R}^{n-1}$ satisfying
\begin{eqnarray}\label{E-thm-q-1-1}
(\mathscr{A}_{1\cdots1}\hat{x}^{m-1})_{i}&=&\sum_{(i_2,\ldots, i_m)\in {\cal N}}a_{i+1i_2\cdots i_m}x_{i_2}\cdots x_{i_m}\nonumber\\
&=&\sum_{(i_2,\ldots, i_m)\in [n]}a_{i+1i_2\cdots i_m}x_{i_2}\cdots x_{i_m}=(\mathscr{A}{x}^{m-1})_{i+1}
\end{eqnarray}
for all $i\in [n-1]$; and
\begin{eqnarray}\label{E-thm-q-1-2}
(\mathscr{A}{x}^{m-1})_1&=&\sum_{(i_2,\ldots, i_m)\in [n]}a_{1i_2\cdots i_m}x_{i_2}\cdots x_{i_m}\nonumber\\
&=&\sum_{(i_2,\ldots, i_m)\in [n]}a_{2i_2\cdots i_m}x_{i_2}\cdots x_{i_m}=(\mathscr{A}{x}^{m-1})_2
\end{eqnarray}
since $a_{1i_{2}\cdots i_{m}}=a_{2i_{2}\cdots i_{m}}$ for all $i_{2},\ldots, i_{m}\in [n]$.

Since $\mathscr{A}_{2\cdots2}\in \mathbb{T}_{m,n}$ is a Q-tensor, it follows that TCP($q_{-1}$,$\mathscr{A}_{1\cdots1}$) has a solution, say $\hat{y}:=(\bar{y}_{2},\ldots, \bar{y}_{n})^{T}$. Then, $\hat{y}\in \mathbb{R}^{n-1}_+$ and for any $i\in [n-1]$,
\begin{eqnarray}\label{E-thm-q-1-3}
0\leq (\mathscr{A}_{1\cdots1}\hat{y}^{m-1}+q_{-1})_i=(\mathscr{A}_{1\cdots1}\hat{y}^{m-1})_{i}+q_{i+1},
\end{eqnarray}
\begin{eqnarray}\label{E-thm-q-1-4}
0&=&\hat{y}_{i}(\mathscr{A}_{1\cdots1}\hat{y}^{m-1}+q_{-1})_{i}=\bar{y}_{i+1}(\mathscr{A}_{1\cdots1}\hat{y}^{m-1})_i+\bar{y}_{i+1}(q_{-1})_{i}\nonumber\\
&=&\bar{y}_{i+1}q_{i+1}+ \bar{y}_{i+1}\sum\limits_{(i_{2},\ldots ,i_{m})\in {\cal N}}a_{i+1i_{2}\cdots i_{m}}\bar{y}_{i_{2}}\cdots \bar{y}_{i_{m}}.
\end{eqnarray}
Let $y:=(0, \hat{y})^T$. Then, we have that $y\in \mathbb{R}^n_+$,
\begin{eqnarray*}
(\mathscr{A}y^{m-1}+q)_{1}&=& (\mathscr{A}y^{m-1})_1+q_{1}\\
&=& (\mathscr{A}y^{m-1})_2+q_{1}\qquad (\mbox{\rm by}\; (\ref{E-thm-q-1-2}))\\
&\geq& (\mathscr{A}y^{m-1})_2+q_{2}\qquad (\mbox{\rm by}\; q_1\geq q_2)\\
&=& (\mathscr{A}_{11\cdots1}\hat{y}^{m-1})_1+q_{2}\qquad (\mbox{\rm by}\; (\ref{E-thm-q-1-1}))\\
&\geq&0,\qquad (\mbox{\rm by}\; (\ref{E-thm-q-1-3}))\\
(\mathscr{A}y^{m-1}+q)_{i+1}&=& (\mathscr{A}_{1\cdots1}\hat{y}^{m-1})_i+q_{i+1}\qquad (\mbox{\rm by}\; (\ref{E-thm-q-1-1}))\\
&\geq&0,\;\; \forall i\in [n-1],\qquad (\mbox{\rm by}\; (\ref{E-thm-q-1-3}))
\end{eqnarray*}
and
\begin{eqnarray*}
y_{1}(\mathscr{A}y^{m-1}+q)_{1}&=& 0,\qquad (\mbox{\rm since}\; y_1=0)\\
y_{i+1}(\mathscr{A}y^{m-1}+q)_{i+1}&=&\bar{y}_{i+1}q_{i+1}+\bar{y}_{i+1}(\mathscr{A}_{1\cdots1}\hat{y}^{m-1})_i\qquad (\mbox{\rm by}\; (\ref{E-thm-q-1-1}))\\
&=&\bar{y}_{i+1}q_{i+1}+\bar{y}_{i+1}\sum\limits_{(i_{2},\ldots ,i_{m})\in {\cal N}}a_{i+1i_{2}\cdots i_{m}}\bar{y}_{i_{2}}\cdots \bar{y}_{i_{m}}\\
&=&0,\;\; \forall i\in [n-1],\qquad (\mbox{\rm by}\; (\ref{E-thm-q-1-4}))
\end{eqnarray*}
Thus, $y$ is a solution to TCP($q,\mathscr{A}$).

Case 2. $q_{1}\leq q_{2}$. In this case, by using the condition that $\mathscr{A}_{1\cdots1}\in \mathbb{T}_{m,n}$ is a Q-tensor, similar to the proof of Case 1, we can obtain that TCP($q,\mathscr{A}$) has a solution.

Combining Case 1 with Case 2, we complete the proof of this theorem.
\ep

\section{Equivalent classes of $Q$-tensors within the class of strong $P_0$-tensors}
\setcounter{equation}{0} \setcounter{Assumption}{0}
\setcounter{Theorem}{0} \setcounter{Proposition}{0}
\setcounter{Corollary}{0} \setcounter{Lemma}{0}
\setcounter{Definition}{0} \setcounter{Remark}{0}
\setcounter{Algorithm}{0}  \setcounter{Example}{0}

\hspace{4mm} In 1979, Agangic and Cottle \cite{ac-79} obtained the following results.
\begin{Proposition}\label{pro-ac}
If $A\in \mathbb{R}^{n\times n}\bigcap P_0$, then
$$ A\in R_0\quad \Longleftrightarrow \quad A\in R\quad \Longleftrightarrow \quad A\in Q.$$
\end{Proposition}
Such a proposition gives two equivalent classes of $Q$-matrices within the class of $P_0$-matrices. After such a pioneer work, the problem on equivalent class of $Q$-matrices has been extensively studied in the literature. See, for example, Pang \cite{Pang-79}, Jeter and Pye \cite{jp-89}, Gowda \cite{Gowda-90}, etc.

In this section, we try to extend the above results to the case of tensor.

Note that both $P_0$-tensor and $P_0^\prime$-tensor are extensions of $P_0$-matrix. It is natural to consider whether the following result holds or not:

{\bf R1}: If $\mathscr{A}\in \mathbb{T}_{m,n}\bigcap P_0$ (or $\mathscr{A}\in \mathbb{T}_{m,n}\bigcap P_0^\prime$), then
 \begin{eqnarray}\label{E-question-1}
 \mathscr{A}\in R_0\quad \Longleftrightarrow \quad \mathscr{A}\in R\quad \Longleftrightarrow \quad \mathscr{A}\in Q.
 \end{eqnarray}
It is regret that the result {\bf R1} does not hold, which can be seen by the following examples.

\begin{Example}\label{exam-p0-1}
Let $\mathscr{A}=(a_{i_1i_2i_3i_4})\in \mathbb{T}_{4,2}$, where $a_{1122}=a_{2222}=1, a_{2112}=-1$ and all other $a_{i_{1}i_{2}i_{3}i_4}=0$. Then $\mathscr{A}\in P_0\bigcap Q$, but $\mathscr{A}\not\in P_0\bigcap R_0$.
\end{Example}

We show that the results in Example \ref{exam-p0-1} hold. Obviously, for any $x\in \mathbb{R}^2$,
$$
\mathscr{A}x^3=\left(\begin{array}{c} x_1x_2^2 \\ x_2^3-x_1^2x_2 \end{array}\right),
$$
and hence,
\begin{eqnarray}\label{E-exam-p0-1}
x_1(\mathscr{A}x^3)_1=x_1^2x_2^2\quad \mbox{\rm and}\quad x_2(\mathscr{A}x^3)_2=x_2^4-x_1^2x_2^2.
\end{eqnarray}
It is easy to see from (\ref{E-exam-p0-1}) that for any $x\in \mathbb{R}^2$ with $x\neq 0$, if $x_1\neq 0$, then $x_1(\mathscr{A}x^3)_1=x_1^2x_2^2\geq 0$; and if $x_1=0$, then $x_2\neq 0$ since $x\neq 0$, and $x_2(\mathscr{A}x^3)_2=x_2^4\geq 0$. Thus, $\mathscr{A}\in P_0$. It is also easy to see from (\ref{E-exam-p0-1}) that $(1,0)^T\in \mathbb{R}^2$ is a solution to TCP$(0,\mathscr{A})$, which, together with the definition of $R_0$-tensor, implies that $\mathscr{A}\not\in R_0$. In the following, we show that $\mathscr{A}\in Q$. Let $a$ and $b$ be two nonnegative real numbers, we consider the following four cases:
\begin{itemize}
\item[C1.] Let $q=(a^{3}, b^{3})^T$, then $z=(0,0)^T$ is a solution to TCP$(q,\mathscr{A})$.
\item[C2.] Let $q=(a^{3}, -b^{3})^T$, then $z=(0, b)^T$ is a solution to TCP$(q,\mathscr{A})$.
\item[C3.] Let $q=(-a^{3}, -b^{3})^T$ with $(a,b)\neq (0,0)$, we show that TCP$(q,\mathscr{A})$ has a solution. In this case, in order to ensure that $(\mathscr{A}x^3)_i+q_i\geq 0$ for $i\in \{1,2\}$, it must hold that $x_{1}\neq 0$ and $x_{2}\neq 0$. So we need to show that the system of equations
    \begin{eqnarray}\label{E-exam-p0-2}
    0=\mathscr{A}x^{3}+q=\left(\begin{array}{c} x_1x_2^2-a^3 \\ x_2^3-x_1^2x_2-b^3 \end{array}\right)
    \end{eqnarray}
    has a nonnegative solution. From the first equation in (\ref{E-exam-p0-2}) it follows that $x_{1}=\frac{a^{3}}{x_{2}^{2}}$; and hence, the second equation in (\ref{E-exam-p0-2}) becomes
    $$(x_{2}^{3})^{2}-b^{3}x_{2}^{3}-a^{6}=0.$$
    It is easy to see that the above equation has a solution $x_{2}^*:=[(b^3+\sqrt{b^6+4a^6})/2]^{1/3}>0$.  Furthermore, $({a^{3}}/({x_{2}^*)^{2}}, x_{2}^*)^T$ is a solution to TCP$(q,\mathscr{A})$.
\item[C4.] Let $q=(-a^{3}, b^{3})^T$. Similar to the proof given in the case C3, we can obtain that TCP$(q,\mathscr{A})$ has a solution in this case.
\end{itemize}
Combining the above four cases, we obtain that $\mathscr{A}\in Q$.

Example \ref{exam-p0-1} demonstrates that (\ref{E-question-1}) cannot be obtained under the assumption that $\mathscr{A}\in \mathbb{T}_{m,n}\bigcap P_0$.

\begin{Example}\label{exam-p0-2}
Let $\mathscr{A}=(a_{i_1i_2i_3})\in \mathbb{T}_{3,2}$, where $a_{122}=a_{222}=1, a_{212}=-1$ and all other $a_{i_{1}i_{2}i_{3}}=0$. Then $\mathscr{A}\in P_0^\prime\bigcap Q$, but $\mathscr{A}\not\in P_0^\prime\bigcap R_0$.
\end{Example}

We show that the results in Example \ref{exam-p0-2} hold.

First, it is obvious that TCP$(0,\mathscr{A})$ is to find $x\in \mathbb{R}^2$ such that
$$
x\geq0,\quad \mathscr{A}x^{2}=\left(\begin{array}{c} x_{2}^{2}\\ x_{2}^{2}-x_{1}x_{2} \end{array}\right)\geq0,\quad x^{T}\mathscr{A}x^{2}=0.
$$
It is easy to see that $(1,0)^T$ is a solution to TCP$(0,\mathscr{A})$; and hence, $\mathscr{A}\not\in R_0$.

Second, we show that $\mathscr{A}$ is a $P_0^\prime$-tensor. For any $x\in \mathbb{R}^2$ with $x\neq0$,
\begin{itemize}
\item if $x_{1}\neq0$, then $x_1^2(\mathscr{A}x^{2})_{1}=x_1^2x_{2}^{2}\geq0$; and
\item if $x_{1}=0$, then $x_2\neq 0$ and $x_2^2(\mathscr{A}x^{2})_{2}=x_{2}^{4}>0$.
\end{itemize}
Thus, $\mathscr{A}$ is a $P_0^\prime$-tensor.

Third, we prove that $\mathscr{A}\in Q$. Let $a$ and $b$ be two nonnegative real numbers.
\begin{itemize}
\item[C1.] Let $q=(a^{2}, b^{2})^T$. Obviously, $(0,0)^T$ is a solution to TCP$(q,\mathscr{A})$.
\item[C2.] Let $q=(-a^{2}, b^{2})^T$ with $a\neq 0$. Take $z:=(({a^{2}+b^{2}})/{a}, a)^T$, then
$$
z\geq0,\;\; \mathscr{A}z^{2}+q=\left(\begin{array}{c} a^{2}-a^{2} \\ a^{2}-a\times \frac{a^{2}+b^{2}}{a}+b^{2}\end{array}\right)=0,\;\; z^T(\mathscr{A}z^{2}+q)=0.
$$
Thus, $z$ solves TCP$(q,\mathscr{A})$ in this case.
\item[C3.] Let $q=(a^{2}, -b^{2})^T$. Take $z:=(0, b)^T$, then
$$
z\geq0,\;\; \mathscr{A}z^{2}+q=\left(\begin{array}{c} b^{2}+a^{2} \\ 0\end{array}\right)\geq0,\;\; z^T(\mathscr{A}z^{2}+q)=0.
$$
Thus, $z$ solves TCP$(q,\mathscr{A})$ in this case.
\item[C4.] Let $q=(-a^{2}, -b^{2})^T$ with $a\leq b$. Take $z:=(0, b)^T$, then
$$
z\geq0,\;\; \mathscr{A}z^{2}+q=\left(\begin{array}{c} b^{2}-a^{2} \\ 0\end{array}\right)\geq0,\;\; z^T(\mathscr{A}z^{2}+q)=0.
$$
Thus, $z$ solves TCP$(q,\mathscr{A})$ in this case.
\item[C5.] Let $q=(-a^{2}, -b^{2})^T$ with $0\neq a\geq b$. Take $z:=(({a^{2}-b^{2}})/{a},a)^T$, then
$$
z\geq0,\;\; \mathscr{A}z^{2}+q=\left(\begin{array}{c} a^{2}-a^{2} \\ a^{2}-a\times\frac{a^{2}-b^{2}}{a}-b^{2}\end{array}\right)=0,\;\; z^T(\mathscr{A}z^{2}+q)=0.
$$
Thus, $z$ solves TCP$(q,\mathscr{A})$ in this case.
\end{itemize}
Therefore, it follows from C1-C5 that TCP$(q,\mathscr{A})$ has a solution for each $q\in \mathbb{R}^2$. Thus, $\mathscr{A}\in Q$.

Example \ref{exam-p0-2} demonstrates that (\ref{E-question-1}) cannot be obtained under the assumption that $\mathscr{A}\in \mathbb{T}_{m,n}\bigcap P_0^\prime$.

Combining Example \ref{exam-p0-1} with Example \ref{exam-p0-2}, we obtain that (\ref{E-question-1}) cannot be obtained within the class of $P_0$-tensors (or $P_0^\prime$-tensors), i.e.,
the result {\bf R1} does not hold.

In order to extend the result of Proposition \ref{pro-ac} to the case of tensor, we introduce a new class of tensors in the following. Recall that the function $f: \mathbb{R}^n\rightarrow \mathbb{R}^n$ is called a $P_0$-function if, for all $x,y\in \mathbb{R}^n$ with $x\neq y$, there is an index $i\in [n]$ such that
$$
x_i\neq y_i\quad\mbox{\rm and}\quad (x_i-y_i)[f_i(x)-f_i(y)]\geq 0.
$$
It is well known that an affine mapping $f(x):=Ax+q$ with $q\in \mathbb{R}^n$ is a $P_0$-function if and only if $A\in \mathbb{R}^{n\times n}$ is a $P_0$-matrix. Inspired by such a result, we introduce a class of tensors which is defined as follows.
\begin{Definition}\label{sp0def}
Given $\mathscr{A}\in \mathbb{T}_{m,n}$. If the mapping $f(x):=\mathscr{A}x^{m-1}+q$ with $q\in \mathbb{R}^n$ is a $P_0$-function, we call $\mathscr{A}$ is a strong $P_0$-tensor, abbreviated as $SP_0$-tensor, and denote the set of all real $m$-order $n$-dimensional $SP_0$-tensors by $\mathbb{T}_{m,n}\bigcap SP_0$.
\end{Definition}

Obviously, when $m=2$, an $SP_0$-tensor reduces to a $P_0$-matrix. Thus, $SP_0$-tensor is an extension of $P_0$-matrix from the matrix space to the tensor space. It is easy to see from Definition \ref{p0def}(i) and Definition \ref{sp0def} that $ \mathbb{T}_{m,n}\bigcap SP_0\subseteq \mathbb{T}_{m,n}\bigcap P_0$.

In the following, we extend the results of Proposition \ref{pro-ac} to the tensor space.
\begin{Theorem}\label{thm-sp0-main}
If $\mathscr{A}\in \mathbb{T}_{m,n}\bigcap SP_0$, we have
 \begin{eqnarray}\label{E-thm-sp0-main-1}
 \mathscr{A}\in R_0\quad \Longleftrightarrow \quad \mathscr{A}\in R \quad \Longleftrightarrow \quad \mathscr{A}\in ER\quad \Longleftrightarrow \quad \mathscr{A}\in Q.
 \end{eqnarray}
\end{Theorem}
{\bf Proof.} Since an $SP_0$-tensor is a $P_0$-tensor and every $P_0$-tensor is semi-positive, it follows from Proposition \ref{pro1}(ii) that
\begin{eqnarray}\label{E-thm-sp0-main-2}
  \mathscr{A}\in R_0 \quad \Longleftrightarrow \quad \mathscr{A}\in R \quad \Longleftrightarrow \quad \mathscr{A}\in ER.
\end{eqnarray}
Thus, in order to show that (\ref{E-thm-sp0-main-1}) holds, we only need to show that
\begin{eqnarray}\label{E-thm-sp0-main-3}
  \mathscr{A}\in R_0 \quad \Longleftrightarrow \quad \mathscr{A}\in Q.
\end{eqnarray}
Suppose that $\mathscr{A}\in R_0$, then $\mathscr{A}\in R$ by (\ref{E-thm-sp0-main-2}). This, together with Proposition \ref{pro1}(i), implies that $\mathscr{A}\in Q$. Thus, in order to show that (\ref{E-thm-sp0-main-3}) holds, we only need to show that $\mathscr{A}\in R_0$ under the condition that $\mathscr{A}\in Q\bigcap SP_0$. Suppose that $\mathscr{A}\in Q$ but $\mathscr{A}\notin R_{0}$. Then there exists a vector $\bar{x}\in \mathbb{R}^n_+\setminus\{0\}$ such that
\begin{eqnarray*}
\left\{\begin{array}{rcl}
(\mathscr{A}\bar{x}^{m-1})_{i}&=&0,\quad\mbox{\rm if}\; \bar{x}_{i}>0, \\
(\mathscr{A}\bar{x}^{m-1})_{i}&\geq&0,\quad\mbox{\rm if}\; \bar{x}_{i}=0.
\end{array}\right.
\end{eqnarray*}
Denote ${\cal I}=\{i\in [n]: \bar{x}_{i}=0\}$ and ${\cal J}=\{i\in [n]: \bar{x}_{i}>0\}$. Take $q\in \mathbb{R}^{n}$  satisfying $q_i>0$ for any $i\in {\cal I}$ and $q_i<0$ for any $i\in {\cal J}$. Since $\mathscr{A}\in Q$, we can assume that $\bar{y}$ is a solution of TCP$(q, \mathscr{A}$). It is obvious that $\bar{x}\neq \bar{y}$. Let $\lambda$ be a positive real number.

For any $i\in \{i\in [n]: \bar{x}_i\neq \bar{y}_i\}$, we consider the following two cases.
\begin{itemize}
\item If $i\in {\cal I}$, then $\bar{x}_{i}=0$ and $\bar{y}_{i}>0$, and hence, it follows that $(\mathscr{A}\bar{x}^{m-1})_{i}\geq0$ and $(\mathscr{A}\bar{y}^{m-1}+q)_{i}=0$. The above equality implies that $(\mathscr{A}\bar{y}^{m-1})_{i}<0$ since $q_{i}>0$ for any $i\in {\cal I}$. This further yields that $(\mathscr{A}(\lambda \bar{y})^{m-1})_{i}=\lambda^{m-1}(\mathscr{A}\bar{y}^{m-1})_{i}<0$ for any $i\in {\cal I}$ since $\lambda>0$. So $(\bar{x}_{i}-\lambda \bar{y}_{i})[(\mathscr{A}\bar{x}^{m-1})_{i}-\mathscr{A}(\lambda \bar{y})^{m-1}_{i}]<0$.
\item If $i\in {\cal J}$, then $(\mathscr{A}\bar{x}^{m-1})_{i}=0$ and $(\mathscr{A}\bar{y}^{m-1}+q)_{i}\geq0$. The above inequality implies that $(\mathscr{A}\bar{y}^{m-1})_{i}>0$ since $q_{i}<0$ for any $i\in {\cal J}$, which yields that $(\mathscr{A}(\lambda \bar{y})^{m-1})_{i}=\lambda^{m-1}(\mathscr{A}\bar{y}^{m-1})_{i}>0$ for any $\lambda>0$. Now, we can choose sufficiently small $\lambda>0$ such that $(\bar{x}-\lambda \bar{y})_{i}>0$. So $(\bar{x}_{i}-\lambda \bar{y}_{i})[(\mathscr{A}\bar{x}^{m-1})_{i}-(\mathscr{A}(\lambda \bar{y})^{m-1})_{i}]<0$ holds for any sufficiently small $\lambda>0$.
\end{itemize}
Thus, we can choose sufficiently small $\lambda>0$ such that for any $i\in \{i\in [n]: \bar{x}_i\neq \bar{y}_i\}$,
$$
(\bar{x}_{i}-\lambda \bar{y}_{i})[(\mathscr{A}\bar{x}^{m-1})_{i}-(\mathscr{A}(\lambda \bar{y})^{m-1})_{i}]<0,
$$
which contradicts the condition that $\mathscr{A}\in SP_{0}$. Therefore, $\mathscr{A}\in R_{0}$; and the desired results are obtained. \ep

In the following, we discuss the relationships among three classes of $P_0$-type tensors.

First, we construct the following example.
\begin{Example}\label{exam-sp0-p0-np0-1}
Let $\mathscr{A}=(a_{i_1\cdots i_m})\in \mathbb{T}_{m,n}$, where $a_{122\cdots 2}=1$ and all other $a_{i_{1}\cdots i_m}=0$. Then $\mathscr{A}\in SP_0\bigcap P_0\bigcap P_0^\prime$.
\end{Example}

We show that the results in Example \ref{exam-sp0-p0-np0-1} hold. Obviously, for any $x\in \mathbb{R}^n$, we have
$$
\mathscr{A}x^{m-1}=\left(\begin{array}{c} x_2^{m-1} \\ 0 \\ \vdots \\ 0 \end{array}\right)\in \mathbb{R}^n.
$$
On one hand, for any $x, y\in \mathbb{R}^n$ with $x\neq y$, if there exists an index $i_0\in \{2,\ldots,n\}$ such that $x_{i_0}\neq y_{i_0}$, then $(x_{i_0}-y_{i_0})[(\mathscr{A}x^{m-1})_{i_0}-(\mathscr{A}y^{m-1})_{i_0}]=0$. Otherwise, we have that $x_i=y_i$ for all $i\in \{2,\ldots,n\}$, and hence, $x_1\neq y_1$ and $x_2=y_2$. Furthermore, $(x_1-y_1)[(\mathscr{A}x^{m-1})_1-(\mathscr{A}y^{m-1})_1]=(x_1-y_1)(x_2^{m-1}-y_2^{m-1})=0$. So $\mathscr{A}\in SP_0$.

On the other hand, for any $x\in \mathbb{R}^n$ with $x\neq 0$, if there exists an index $i_0\in \{2,\ldots,n\}$ such that $x_{i_0}\neq 0$, then $x_{i_0}(\mathscr{A}x^{m-1})_{i_0}=0$ and $x_{i_0}^{m-1}(\mathscr{A}x^{m-1})_{i_0}=0$. Otherwise, we have that $x_i=0$ for all $i\in \{2,\ldots,n\}$, and hence, $x_1\neq 0$ and $x_2=0$. Furthermore, $x_1(\mathscr{A}x^{m-1})_1=x_1x_2^{m-1}=0$ and $x_1^{m-1}(\mathscr{A}x^{m-1})_1=x_1^{m-1}x_2^{m-1}=0$. So $\mathscr{A}\in P_0\bigcap P_0^\prime$.

\begin{Proposition}\label{pro-sp0-p0-np0}
We have that $\mathbb{T}_{m,n}\bigcap SP_0\bigcap P_0\bigcap P_0^\prime\neq \emptyset$.
\end{Proposition}
{\bf Proof.} The results of the proposition hold directly from Example \ref{exam-sp0-p0-np0-1}. \ep

Second, we consider the relationship between the class of $P_0$-tensors and the class of $SP_0$-tensor. We construct the following example.
\begin{Example}\label{exam-p0-sp0-1}
Let $\mathscr{A}=(a_{i_1i_2i_3i_4})\in \mathbb{T}_{4,2}$, where $a_{1122}=a_{2122}=1$ and all other $a_{i_{1}i_{2}i_{3}i_4}=0$. Then $\mathscr{A}\in P_0$, but $\mathscr{A}\not\in SP_0$.
\end{Example}

We show that the results in Example \ref{exam-p0-sp0-1} hold. Obviously, for any $x\in \mathbb{R}^2$ with $x\neq 0$, we have
$$
\mathscr{A}x^3=\left(\begin{array}{c} x_1x_2^2 \\ x_1x_2^2 \end{array}\right).
$$
If $x_{1}\neq0$, then $x_{1}(\mathscr{A}x^{3})_{1}=x_{1}^{2}x_{2}^{2}\geq0$; and if $x_{1}=0$, then $x_{2}\neq0$, and $x_{2}(\mathscr{A}x^{3})_{2}=x_1x_2^3=0$. So we obtain that $\mathscr{A}\in P_0$.

In addition, for any given $q\in \mathbb{R}^2$, let $f(x)=\mathscr{A}x^3+q$ for any $x\in \mathbb{R}^2$; and
take $\bar{x}=(1,1)^T$ and $\bar{y}=(1,-2)^T$, then it is easy to see that
$$
\bar{x}_{1}=\bar{y}_{1},\quad \mbox{\rm and}\quad \bar{x}_{2}\neq \bar{y}_{2},\; (\bar{x}_{2}-\bar{y}_{2})(f(\bar{x})-f(\bar{y}))_{2}=-9<0.
$$
These demonstrate that $\mathscr{A}\not\in SP_0$.

\begin{Proposition}\label{pro-p0-sp0}
If $\mathscr{A}$ is an $SP_0$-tensor, then it is a $P_0$-tensor. But the converse is not true.
\end{Proposition}
{\bf Proof.} The first result holds from the definition of $SP_0$-tensor given in Definition \ref{sp0def} and the definition of $P_0$-tensor given in Definition \ref{p0def}(i); and the second result holds from Example \ref{exam-p0-sp0-1}. \ep

Third, we consider the relationship between the class of $P_0$-tensors and the class of $P_0^\prime$-tensors. We construct the following two examples.

\begin{Example}\label{exam-p0-np0-1}
Let $\mathscr{A}=(a_{i_1i_2i_3})\in \mathbb{T}_{3,2}$, where $a_{121}=1, a_{211}=-1$ and all other $a_{i_{1}i_{2}i_{3}}=0$. Then $\mathscr{A}\in P_0$, but $\mathscr{A}\not\in P_0^\prime$.
\end{Example}

We show that the results in Example \ref{exam-p0-np0-1} hold. Obviously, for any $x\in \mathbb{R}^2$, we have
$$
\mathscr{A}x^2=\left(\begin{array}{c} x_2x_1 \\ -x_1^2 \end{array}\right).
$$
On one hand, from $x_1(\mathscr{A}x^2)_1=x_2x_1^2$ and $x_2(\mathscr{A}x^2)_2=-x_2x_1^2$,
it is easy to see that for any $x\in \mathbb{R}^2$ with $x\neq 0$, there exists an index $i\in \{1,2\}$ such that $x_i\neq 0$ and $x_i(\mathscr{A}x^{2})_i\geq 0$, i.e., $\mathscr{A}\in P_0$. On the other hand, for any $\alpha>0$ and $\beta<0$, by taking $(x_1,x_2):=(\alpha,\beta)$, we have
$$
x_1^2(\mathscr{A}x^2)_1=x_2x_1^3=\beta\alpha^3<0\quad \mbox{\rm and}\quad x_2^2(\mathscr{A}x^2)_2=-x_2^2x_1^2=-\beta^2\alpha^2<0,
$$
and hence, $\mathscr{A}\not\in P_0^\prime$.

\begin{Example}\label{exam-p0-np0-2}
Let $\mathscr{A}=(a_{i_1i_2i_3})\in \mathbb{T}_{3,2}$, where $a_{122}=1, a_{211}=-1$ and all other $a_{i_{1}i_{2}i_{3}}=0$. Then $\mathscr{A}\in P_0^\prime$, but $\mathscr{A}\not\in P_0$.
\end{Example}

We show that the results in Example \ref{exam-p0-np0-2} hold. Obviously, for any $x\in \mathbb{R}^2$, we have
$$
\mathscr{A}x^2=\left(\begin{array}{c} x_2^2 \\ -x_1^2 \end{array}\right).
$$
On one hand, from $x_1^2(\mathscr{A}x^2)_1=x_1^2x_2^2$ and $x_2^2(\mathscr{A}x^2)_2=-x_2^2x_1^2$,
it is easy to see that for any $x\in \mathbb{R}^2$ with $x\neq 0$, there exists an index $i\in \{1,2\}$ such that $x_i\neq 0$ and $x_i^2(\mathscr{A}x^{2})_i\geq 0$, i.e., $\mathscr{A}\in P_0^\prime$. On the other hand, for any $\alpha<0$ and $\beta>0$, by taking $(x_1,x_2):=(\alpha,\beta)$, we have
$$
x_1(\mathscr{A}x^2)_1=x_1x_2^2=\alpha \beta^2<0\quad \mbox{\rm and}\quad x_2(\mathscr{A}x^2)_2=-2x_2x_1^2=-\beta\alpha^2<0,
$$
and hence, $\mathscr{A}\not\in P_0$.

\begin{Proposition}\label{pro-p0-np0}
(i) We have $\mathbb{T}_{m,n}\bigcap P_0^\prime=\mathbb{T}_{m,n}\bigcap P_0$ when $m$ is even. (ii) There is no inclusion relation between the class of odd order $P_0$-tensors and the class of odd order $P_0^\prime$-tensors.
\end{Proposition}
{\bf Proof.} The result (i) holds directly from the definitions of $P_0$-tensor and $P_0^\prime$-tensor given in Definition \ref{p0def}; and the result (ii) holds directly from Examples \ref{exam-p0-np0-1} and \ref{exam-p0-np0-2}. \ep

Fourth, we consider the relationship between the class of $SP_0$-tensor and the class of $P_0^\prime$-tensor. From Proposition \ref{pro-p0-sp0} and Proposition \ref{pro-p0-np0}(i), we obtain immediately the following results.

\begin{Proposition}\label{pro-sp0-np0}
When $m$ is even, every $m$-order $SP_0$-tensor is an $m$-order $P_0^\prime$-tensor, but the converse is not true.
\end{Proposition}

In the following, we consider the odd order $SP_0$-tensors and odd order $P_0^\prime$-tensors.

\begin{Lemma}\label{Lem-sp0-1}
Let $\mathscr{A}\in \mathbb{T}_{m,n}\bigcap SP_0$ with $m$ being odd. Then, for any $i\in [n]$, we have that either $(\mathscr{A}x^{m-1})_i\equiv 0$ or $(\mathscr{A}x^{m-1})_i$ is a function of variables $x_1,\ldots, x_{i-1}, x_{i+1}, \ldots, x_n$, but independent of the variable $x_i$.
\end{Lemma}
{\bf Proof.} For any $x=(x_1, \ldots, x_n)^T\in \mathbb{R}^n$, suppose that $a$ is an arbitrary fixed real number and $i_0\in [n]$, we take $y=(x_1, \ldots, x_{i_0-1}, a, x_{i_0+1}, \ldots, x_n)^T$, then for any
$$
i\in {\cal I}:=\{1,\ldots,i_0-1,i_0+1,\ldots,n\},
$$
we have $x_i=y_i$. For any $x_{i_0}\in \mathbb{R}\setminus\{a\}$, we have $x_{i_0}\neq y_{i_0}$, which, together with $\mathscr{A}\in SP_0$, implies that
\begin{eqnarray}\label{E-Lem-sp0-1-1}
(x_{i_0}-a)[(\mathscr{A}x^{m-1})_{i_0}-(\mathscr{A}y^{m-1})_{i_0}]\geq0.
\end{eqnarray}
For $-x$ and $-y$, we have $\mathscr{A}(-x)^{m-1}=\mathscr{A}x^{m-1}$ since $m$ is odd; and $-x_{i_0}\neq -y_{i_0}$ and $-x_{i}=-y_{i}$ for any $i\in {\cal I}$. These and $\mathscr{A}\in SP_0$ imply that
\begin{eqnarray}\label{E-Lem-sp0-1-2}
(a-x_{i_0})[(\mathscr{A}x^{m-1})_{i_0}-(\mathscr{A}y^{m-1})_{i_0}]\geq0.
\end{eqnarray}
Combining (\ref{E-Lem-sp0-1-1}) with (\ref{E-Lem-sp0-1-2}), we obtain that for any $x_{i_0}\in \mathbb{R}\setminus\{a\}$,
$$
(\mathscr{A}x^{m-1})_{i_0}=(\mathscr{A}y^{m-1})_{i_0}.
$$
By the arbitrariness of $x_{i_0}$, the above equality implies that either $(\mathscr{A}x^{m-1})_{i_0}\equiv 0$ or $(\mathscr{A}x^{m-1})_{i_0}$ is independent of the variable $x_{i_0}$. Furthermore, the desired results holds by the arbitrariness of ${i_0}$. \ep

\begin{Proposition}\label{pro-sp0-2}
Let $m$ be odd. Then $\mathbb{T}_{m,2}\bigcap SP_0\subseteq \mathbb{T}_{m,2}\bigcap P_0^\prime$.
\end{Proposition}
{\bf Proof.} Given $\mathscr{A}\in \mathbb{T}_{m,2}\bigcap SP_0$. Then $\mathscr{A}\in \mathbb{T}_{m,2}\bigcap P_0$.

First, we show that for any $x\in \mathbb{R}^2$, there exists an index $i\in \{1,2\}$ such that $(\mathscr{A}x^{m-1})_i\equiv 0$. We assume that $(\mathscr{A}x^{m-1})_i\not\equiv 0$ for all $i\in \{1,2\}$. Since $\mathscr{A}\in \mathbb{T}_{m,2}\bigcap SP_0$, it follows from Lemma \ref{Lem-sp0-1} that
$$
\mathscr{A}x^{m-1}=\left(\begin{array}{c} \alpha x_2^{m-1} \\ \beta x_1^{m-1}  \end{array}\right),
$$
where $\alpha,\beta\in \mathbb{R}\setminus\{0\}$. Without loss of generality, we assume that $\bar{x}=(\bar{x}_1, \bar{x}_2)^T\in \mathbb{R}^2_+$ such that $(\mathscr{A}x^{m-1})_i\neq 0$ for all $i\in \{1,2\}$. Take
\begin{eqnarray*}
\hat{x}=\left\{\begin{array}{ll}
(-\bar{x}_1, -\bar{x}_2)^T\quad &\mbox{\rm if}\; (\mathscr{A}\bar{x}^{m-1})_1>0\; \mbox{\rm and}\; (\mathscr{A}\bar{x}^{m-1})_2>0,\\
(-\bar{x}_1, \bar{x}_2)^T\quad &\mbox{\rm if}\; (\mathscr{A}\bar{x}^{m-1})_1>0\; \mbox{\rm and}\; (\mathscr{A}\bar{x}^{m-1})_2<0,\\
(\bar{x}_1, -\bar{x}_2)^T\quad &\mbox{\rm if}\; (\mathscr{A}\bar{x}^{m-1})_1<0\; \mbox{\rm and}\; (\mathscr{A}\bar{x}^{m-1})_2>0,\\
(\bar{x}_1, \bar{x}_2)^T\quad &\mbox{\rm if}\; (\mathscr{A}\bar{x}^{m-1})_1<0\; \mbox{\rm and}\; (\mathscr{A}\bar{x}^{m-1})_2<0,\\
\end{array}\right.
\end{eqnarray*}
then $\mathscr{A}\hat{x}^{m-1}=\mathscr{A}\bar{x}^{m-1}$ since $m$ is odd; and hence, $\hat{x}_1(\mathscr{A}\hat{x}^{m-1})_1<0$ and $\hat{x}_2(\mathscr{A}\hat{x}^{m-1})_2<0$, which is a contradiction with $\mathscr{A}\in P_0$.

Second, we show that $\mathscr{A}\in \mathbb{T}_{m,2}\bigcap P_0^\prime$. Without loss of generality, we assume that $(\mathscr{A}x^{m-1})_2\equiv0$ for any $x\in \mathbb{R}^2$.  For any $x\in \mathbb{R}^2$ with $x\neq 0$,
\begin{itemize}
\item if $x_2\neq0$, then $x_2^{m-1}(\mathscr{A}x^{m-1})_2=0$; and
\item if $x_2=0$, then $x_1\neq 0$ and $x_1^{m-1}(\mathscr{A}x^{m-1})_1=x_1^{m-1}(\alpha x_2^{m-1})=0$,
\end{itemize}
so $\mathscr{A}\in P_0^\prime$. Therefore, the desired result holds. \ep

From Example \ref{exam-sp0-p0-np0-1} and Proposition \ref{pro-sp0-2}, we have obtained some relationship between $\mathbb{T}_{m,n}\bigcap SP_0$ and $\mathbb{T}_{m,n}\bigcap P_0^\prime$. But this does not give a full characterization for the relationship between $\mathbb{T}_{m,n}\bigcap SP_0$ and $\mathbb{T}_{m,n}\bigcap P_0^\prime$. We conjecture that it is possible that the class of $SP_0$-tensors is a proper subset of the class of $P_0^\prime$-tensors, which needs to be further studied in the future.

At the end of this section, we give a characterization of $SP_0$-tensor. In \cite{sq-15}, Song and Qi gave the concept of principal sub-tensors and proved that the principal sub-tensors of every $P_0$-tensor is a $P_0$-tensor. A tensor $\mathscr{C}\in \mathbb{T}_{m,r}$ is called a principal sub-tensor of tensor $\mathscr{A}=(a_{i_1i_2\cdots i_m})\in \mathbb{T}_{m,n}\;(1\leq r\leq n)$ if there is a set ${\cal J}$ that composed of $r$ elements in $[n]$ such that $\mathscr{C}=(a_{i_1i_2\cdots i_m})$ for all $i_1, \ldots ,i_m\in {\cal J}$.

\begin{Theorem}\label{thm-sp0-1}
Let $\mathscr{A}\in \mathbb{T}_{m,n}$ be an $SP_0$-tensor, then every principal sub-tensor of $\mathscr{A}$ is an $SP_0$-tensor.
\end{Theorem}
{\bf Proof.} Let $\mathscr{A}_r^{\cal J}$ be an arbitrary principal sub-tensor of $\mathscr{A}$. Suppose that $\mathscr{A}_r^{\cal J}$ is not an $SP_0$-tensor, then for each pair of distinct vectors $x=(x_{j_1},\ldots, x_{j_r})^T\in \mathbb{R}^r$ and $y=(y_{j_1},\ldots, y_{j_r})^T\in \mathbb{R}^r$, it follows that the set ${\cal K}:=\{j\in {\cal J}: x_j\neq y_j\}$ is nonempty and
$$
(x_j-y_j)(\mathscr{A}_r^{\cal J}x^{m-1}-\mathscr{A}_r^{\cal J}y^{m-1})_j<0,\quad \forall j\in {\cal K}.
$$
Let $x^*=(x_1^*,\ldots,x_n^*)^T\in \mathbb{R}^n$ and $y^*=(y_1^*,\ldots,y_n^*)^T\in \mathbb{R}^n$ be defined by
\begin{eqnarray*}
x_i^*=\left\{\begin{array}{ll} x_i\;\; &\mbox{\rm if}\; i\in {\cal J},\\ 0\;\; &\mbox{\rm otherwise}, \end{array}\right. \quad
y_i^*=\left\{\begin{array}{ll} y_i\;\; &\mbox{\rm if}\; i\in {\cal J},\\ 0\;\; &\mbox{\rm otherwise}, \end{array}\right. \quad
\forall i\in [n],
\end{eqnarray*}
then it is easy to show that
$$
(x^*_i-y^*_i)[(\mathscr{A}(x^*)^{m-1})_i-(\mathscr{A}(y^*)^{m-1})_i]<0
$$
holds for any $i\in \{i\in [n]: x^*_i\neq y^*_i\}$. This contradicts the condition that $\mathscr{A}\in SP_0$. Thus, the desired result holds.\ep

\section{Equivalent classes of $Q$-tensors within the class of nonegative tensors}
\setcounter{equation}{0} \setcounter{Assumption}{0}
\setcounter{Theorem}{0} \setcounter{Proposition}{0}
\setcounter{Corollary}{0} \setcounter{Lemma}{0}
\setcounter{Definition}{0} \setcounter{Remark}{0}
\setcounter{Algorithm}{0} \setcounter{Example}{0}

\hspace{4mm} The following result was obtained by Danao in \cite[Theorem 4.8]{Danao-94}.
\begin{Proposition}\label{pro-d}
If $A\in \mathbb{R}^{n\times n}$ is nonegative, then
$$A\in R_0\quad \Longleftrightarrow\quad A\in Q.$$
\end{Proposition}

A natural question is whether the result of Proposition \ref{pro-d} can be extended to the case of tensor or not. In this section, we give a positive answer to this question, which is given as follows.
\begin{Theorem}\label{thm-nn-main}
If $\mathscr{A}\in \mathbb{T}_{m,n}$ is nonnegative, we have
 \begin{eqnarray}\label{E-thm-nn-main-1}
 \mathscr{A}\in R_0\quad \Longleftrightarrow \quad \mathscr{A}\in R\quad \Longleftrightarrow \quad \mathscr{A}\in Q \quad \Longleftrightarrow \quad \mathscr{A}\in ER.
 \end{eqnarray}
\end{Theorem}
{\bf Proof.} We consider the following five cases.

(a) We show that $\mathscr{A}\in R_0\; \Rightarrow \; \mathscr{A}\in R$. Suppose that $\mathscr{A}$ is not an $R$-tensor. By the definition of $R$-tensor, there exists  $(\hat{x},\hat{t})\in (\mathbb{R}^n_+\setminus\{0\})\times \mathbb{R}_+$ such that
\begin{eqnarray*}
\left\{\begin{array}{ll}
(\mathscr{A}\hat{x}^{m-1})_{i}+\hat{t}=0,\quad & \mbox{\rm if}\; \hat{x}_{i}>0,\vspace{2mm}\\
(\mathscr{A}\hat{x}^{m-1})_{i}+\hat{t}\geq 0,\quad & \mbox{\rm if}\; \hat{x}_{i}=0.
\end{array}\right.
\end{eqnarray*}
If $\hat{t}=0$, then the above contradicts the condition that $\mathscr{A}\in R_{0}$. So $\hat{t}>0$. Thus, for any $i\in [n]$, we have
\begin{eqnarray*}
\left\{\begin{array}{ll}
(\mathscr{A}\hat{x}^{m-1})_{i}=-\hat{t}<0,\quad & \mbox{\rm if}\; \hat{x}_{i}>0,\vspace{2mm}\\
(\mathscr{A}\hat{x}^{m-1})_{i}+\hat{t}\geq 0,\quad & \mbox{\rm if}\; \hat{x}_{i}=0,
\end{array}\right.
\end{eqnarray*}
which yields
 \begin{eqnarray}\label{E-thm-nn-main-2}
 \mathscr{A}\hat{x}^m=\hat{x}^T\mathscr{A}\hat{x}^{m-1}
 =\sum_{i=1}^n\hat{x}_i(\mathscr{A}(\hat{x})^{m-1})_i<0.
 \end{eqnarray}
In addition, since $\hat{x}\in \mathbb{R}^n_+\setminus\{0\}$ and $\mathscr{A}$ is a nonnegative tensor, i.e., $a_{i_{1}\cdots i_{m}}\geq 0$ for all $i_{1},\ldots, i_{m}\in [n]$, we have
$$
\mathscr{A}\hat{x}^m=\sum_{i_1,\ldots,i_{m}=1}^{n}a_{i_1\cdots i_{m}}\hat{x}_{i_{1}}\cdots\hat{x}_{i_{m}}\geq 0,
$$
which contradicts (\ref{E-thm-nn-main-2}). So $\mathscr{A}\in R_{0}\Rightarrow \mathscr{A}\in R$.

(b) We show that $\mathscr{A}\in R\; \Rightarrow \; \mathscr{A}\in Q$. Such a result holds from Proposition \ref{pro1}(i).

(c) We show that $\mathscr{A}\in Q\; \Rightarrow \; \mathscr{A}\in R_0$. Suppose that $\mathscr{A}\not\in R_0$. Then, by the definition of $R_0$-tensor, there exists  $\hat{x}\in \mathbb{R}^n_+$ with $\hat{x}\neq 0$ such that
\begin{eqnarray*}
\left\{\begin{array}{ll}
(\mathscr{A}\hat{x}^{m-1})_{i}=0,\quad & \mbox{\rm if}\; \hat{x}_{i}>0,\vspace{2mm}\\
(\mathscr{A}\hat{x}^{m-1})_{i}\geq 0,\quad & \mbox{\rm if}\; \hat{x}_{i}=0.
\end{array}\right.
\end{eqnarray*}
Thus, there exists $j\in [n]$ such that $\hat{x}_{j}>0$, and then, $(\mathscr{A}\hat{x}^{m-1})_{j}=0$. Since $\mathscr{A}\in Q$, it follows from Proposition \ref{pro1}(iv) that $a_{j\cdots j}>0$. So,
$$
0=(\mathscr{A}\hat{x}^{m-1})_{j}=\sum_{i_{2},\ldots,i_{m}=1}^{n}a_{ji_{2}\cdots i_{m}}\hat{x}_{i_{2}}\cdots\hat{x}_{i_{m}}\geq a_{j\cdots j}x_{j}\cdots x_{j}>0,
$$
which derives a contradiction. Thus, $\mathscr{A}\in Q\Rightarrow \mathscr{A}\in R_{0}$.

(d) We show that $\mathscr{A}\in R_0\; \Rightarrow \; \mathscr{A}\in ER$. Suppose that $\mathscr{A}$ is not an $ER$-tensor. By the definition of $ER$-tensor, there exists  $(\hat{x},\hat{t})\in (\mathbb{R}^n_+\setminus\{0\})\times \mathbb{R}_+$ such that
\begin{eqnarray*}
\left\{\begin{array}{lcll}
(\mathscr{A}\hat{x}^{m-1})_{i}+\hat{t}\hat{x}_i&=&0,\quad & \mbox{\rm if}\; \hat{x}_{i}>0,\vspace{2mm}\\
(\mathscr{A}\hat{x}^{m-1})_{i}   &\geq& 0,\quad & \mbox{\rm if}\; \hat{x}_{i}=0.
\end{array}\right.
\end{eqnarray*}
If $\hat{t}=0$, then the above contradicts the condition that $\mathscr{A}\in R_{0}$. So $\hat{t}>0$. Denote ${\cal I}:=\{i\in [n]: \hat{x}_{i}>0\}$, then
 \begin{eqnarray}\label{E-thm-nn-main-3}
 \mathscr{A}\hat{x}^m=\hat{x}^T\mathscr{A}\hat{x}^{m-1}
 =\sum_{i=1}^n\hat{x}_i(\mathscr{A}\hat{x}^{m-1})_i
 =\sum_{i\in {\cal I}}\hat{x}_i(-\hat{t}\hat{x}_i)<0.
 \end{eqnarray}
In addition, since $\hat{x}\in \mathbb{R}^n_+\setminus\{0\}$ and $\mathscr{A}$ is a nonnegative tensor, we have
$$
\mathscr{A}\hat{x}^m=\sum_{i_1,\ldots,i_{m}=1}^{n}a_{i_1\cdots i_{m}}\hat{x}_{i_{1}}\cdots\hat{x}_{i_{m}}\geq 0,
$$
which contradicts (\ref{E-thm-nn-main-3}). So $\mathscr{A}\in R_{0}\Rightarrow \mathscr{A}\in ER$.

(e) We show that $\mathscr{A}\in ER\; \Rightarrow \; \mathscr{A}\in Q$. Such a result holds from Proposition \ref{pro1}(iii).

Combining cases (a)-(e), we obtain that (\ref{E-thm-nn-main-1}) holds within the class of nonnegative tensors. \ep

In Theorems \ref{thm-sp0-main} and \ref{thm-nn-main}, we have obtained several equivalent classes of $Q$-tensors within the class of $SP_0$-tensors and the class of nonnegative tensors, respectively. What is the relationship between the class of $SP_{0}$-tensors and the class of nonnegative tensors? We see the following two examples.

\begin{Example}\label{exam-sp0-nonnegative-1}
Let $\mathscr{A}=(a_{i_1i_2i_3})\in \mathbb{T}_{3,2}$, where $a_{111}=1, a_{222}=1$ and all other $a_{i_{1}i_{2}i_{3}}=0$. Then $\mathscr{A}$ is nonnegative, but $\mathscr{A}\not\in SP_0$.
\end{Example}

We show that the results in Example \ref{exam-sp0-nonnegative-1} hold. On one hand, it is obvious that $\mathscr{A}$ is nonnegative. On the other hand, we have
$$
\mathscr{A}x^2=\left(\begin{array}{c} x_1^2 \\ x_2^2 \end{array}\right),\quad \forall x\in \mathbb{R}^2.
$$
Take $x=(-2,-3)^T$, $y=(1,2)^T$, then
$$
\max\left\{(x_{i}-y_{i})(\mathscr{A}x^2-\mathscr{A}y^2)_i: i\in \{1,2\}\right\}=\max\{-9,-12\}=-9<0,
$$
and hence, $\mathscr{A}\not\in SP_0$.

\begin{Example}\label{exam-sp0-nonnegative-2}
Let $\mathscr{A}=(a_{i_1i_2i_3i_4})\in \mathbb{T}_{4,2}$, where $a_{1122}=-1, a_{2222}=1$ and all other $a_{i_{1}i_{2}i_{3}i_4}=0$. Then $\mathscr{A}$ is not nonnegative, but $\mathscr{A}\in SP_0$.
\end{Example}

We show that the results in Example \ref{exam-sp0-nonnegative-2} hold. On one hand, it is obvious that $\mathscr{A}$ is not nonnegative. On the other hand, we have
$$
\mathscr{A}x^3=\left(\begin{array}{c} -x_1x_2^2 \\ x_2^3 \end{array}\right),\quad \forall x\in \mathbb{R}^2.
$$
For any $x,y\in \mathbb{R}^2$ with $x\neq y$, we have
\begin{eqnarray*}
\max\left\{(x_i-y_i)[(\mathscr{A}x^3)_i-(\mathscr{A}y^3)_i)]: i\in \{1,2\}\right\}&\geq& (x_2-y_2)[(\mathscr{A}x^3)_2-(\mathscr{A}y^3)_2]\\
&=&(x_2-y_2)(x_2^3-y_2^3)\\
&=&(x_2-y_2)^2(x_2^2+x_2y_2+y_2^2)\\
&\geq&0,
\end{eqnarray*}
and hence, $\mathscr{A}\in SP_0$.

From Examples \ref{exam-sp0-nonnegative-1} and \ref{exam-sp0-nonnegative-2}, it follows that the class of nonnegative tensors is different from the class of $SP_0$-tensors. Combining Theorems \ref{thm-sp0-main} and \ref{thm-nn-main}, we have obtained that (\ref{E-thm-nn-main-1}) holds within the class of nonnegative tensors or the class of $SP_0$-tensors.

\section{Q-tensors, semi-positive tensors and copositive tensors}
\setcounter{equation}{0} \setcounter{Assumption}{0}
\setcounter{Theorem}{0} \setcounter{Proposition}{0}
\setcounter{Corollary}{0} \setcounter{Lemma}{0}
\setcounter{Definition}{0} \setcounter{Remark}{0}
\setcounter{Algorithm}{0} \setcounter{Example}{0}

\hspace{4mm} In this section, we show that three known results related to $Q$-matrices cannot be extended to the tensor space by using several examples.

First, we consider the following result which was obtained by Pang \cite{Pang-79}.
\begin{Proposition}\label{pro-p}
Let $A\in \mathbb{R}^{n\times n}\bigcap L_1\bigcap Q$. If $x^*$ is a nonzero solution to LCP$(0,A)$, then $x^*$ contains at least two nonzero components.
\end{Proposition}

Since it is not clear whether such a result can be extended to the case of tensor or not, Song and Qi proposed the following question (see Question 3.1 in \cite{sq-15-r1}):

\noindent {\bf Q}: Whether or not a nonzero solution $x$ of TCP$(0,\mathscr{A})$ contains at least two nonzero components if $\mathscr{A}$ is a semi-positive Q-tensor.

We now answer this question by using Example \ref{exam-p0-2}. Let $\mathscr{A}$ be given by Example \ref{exam-p0-2}. On one hand, it has been showed that $\mathscr{A}\in P_0^\prime\bigcap Q$, which implies that $\mathscr{A}$ is a semi-positive $Q$-tensor since every $P_0^\prime$-tensor is a semi-positive tensor. On the other hand, it is obvious that $(1,0)^T$ is a solution to TCP$(0,\mathscr{A})$.
Therefore, by Example \ref{exam-p0-2} we obtain that the results of Proposition \ref{pro-p} cannot be extended to the case of tensor.

Second, we consider the following result which was obtained by Pang \cite{Pang-79}.
\begin{Proposition}\label{pro-p-2}
Let $A\in \mathbb{R}^{n\times n}\bigcap L_1\bigcap Q$. Then the system
$$
Mx=0,\quad x>0
$$
is inconsistent.
\end{Proposition}

A natural question is whether or not the system
\begin{eqnarray}\label{E-q3-1}
\mathscr{A}x^{m-1}=0,\quad x>0
\end{eqnarray}
is inconsistent for any semi-positive $\mathscr{A}\in \mathbb{T}_{m,n}\bigcap Q$.
We now answer this question by constructing the following example.

\begin{Example}\label{exam-q-l1-2}
Let $m\geq 3$ be odd and $\mathscr{A}=(a_{i_1\cdots i_m})\in \mathbb{T}_{m,2}$, where
$$
a_{\tiny \underbrace{1\cdots 1}_{m-i}\underbrace{2\cdots 2}_i}=a_{\tiny \underbrace{2\cdots 2}_{i+1}\underbrace{1\cdots 1}_{m-1-i}}=(-1)^iC^{i}_{m-1},\quad \forall i\in \{1,\ldots,m-1\}
$$
and all other $a_{i_1\cdots i_m}=0$. Then $\mathscr{A}\in Q$ is semi-positive, but $(1,1)^T$ is a solution to the system (\ref{E-q3-1}).
\end{Example}

We show that the results in Example \ref{exam-q-l1-2} hold. Obviously, for any $x\in \mathbb{R}^n$, it follows that
$$
\mathscr{A}x^{m-1}=\left(\begin{array}{c} (x_1-x_2)^{m-1} \\ (x_1-x_2)^{m-1}\end{array}\right).
$$
It is easy to see that $(1,1)^T$ is a solution to the system (\ref{E-q3-1}). Moreover, for any $x\in \mathbb{R}^2_+$ with $x\neq0$,
\begin{itemize}
\item if $x_{1}>0$, then $x_1(\mathscr{A}x^{2})_{1}=x_1(x_1-x_2)^{m-1}\geq0$; and
\item if $x_{1}=0$, then $x_{2}>0$ and $x_2(\mathscr{A}x^{2})_{2}=x_2(x_1-x_2)^{m-1}>0$.
\end{itemize}
Thus, $\mathscr{A}$ is a semi-positive tensor. In the following, we prove that $\mathscr{A}\in Q$. For any $a,b\in \mathbb{R}_+$, we consider the following five cases.
\begin{itemize}
\item[C1.] Let $q=(a^{m-1}, b^{m-1})^T$. Obviously, $(0,0)^T$ is a solution to TCP$(q,\mathscr{A})$.
\item[C2.] Let $q=(-a^{m-1}, b^{m-1})^T$. Take $z:=(a, 0)^T$, then
$$
z\geq0,\;\; \mathscr{A}z^{m-1}+q=\left(\begin{array}{c} 0 \\ a^{m-1}+b^{m-1}\end{array}\right)\geq 0,\;\;z^T(\mathscr{A}z^{m-1}+q)=0.
$$
Thus, $z$ solves TCP$(q,\mathscr{A})$ in this case.
\item[C3.] Let $q=(a^{m-1}, -b^{m-1})^T$. Take $z:=(0, b)^T$, then
$$
z\geq0,\;\; \mathscr{A}z^{m-1}+q=\left(\begin{array}{c}a^{m-1}+b^{m-1} \\ 0\end{array}\right)\geq0,\;\;z^T(\mathscr{A}z^{m-1}+q)=0.
$$
Thus, $z$ solves TCP$(q,\mathscr{A})$ in this case.
\item[C4.] Let $q=(-a^{m-1}, -b^{m-1})^T$ with $a\leq b$. Take $z:=(0, b)^T$, then
$$
z\geq0,\;\; \mathscr{A}z^{m-1}+q=\left(\begin{array}{c} b^{m-1}-a^{m-1} \\ 0\end{array}\right)\geq0,\;\;z^T(\mathscr{A}z^{m-1}+q)=0.
$$
Thus, $z$ solves TCP$(q,\mathscr{A})$ in this case.
\item[C5.] Let $q=(-a^{m-1}, -b^{m-1})^T$ with $a\geq b$. Take $z:=(a, 0)^T$, then
$$
z\geq0,\;\; \mathscr{A}z^{m-1}+q=\left(\begin{array}{c} 0 \\ a^{m-1}-b^{m-1}\end{array}\right)\geq 0,\;\;z^T(\mathscr{A}z^{m-1}+q)=0.
$$
Thus, $z$ solves TCP$(q,\mathscr{A})$ in this case.
\end{itemize}
Combining cases C1-C5, we obtain that TCP$(q,\mathscr{A})$ has a solution for each $q\in \mathbb{R}^2$. Thus, $\mathscr{A}\in Q$.

Example \ref{exam-q-l1-2} shows that the result of Proposition \ref{pro-p-2} cannot be extended to the case of tensor.

Third, we consider the following result which was obtained by Jeter and Pye \cite{jp-85}.
\begin{Proposition}\label{pro-p-22}
If $M\in \mathbb{R}^{n\times n}$ is copositive and $n\leq3$, then
$$M\in Q\quad \Longleftrightarrow \quad M\in R_0.$$
\end{Proposition}

A natural question is whether the above result can be extended to the tensor space or not. In order to answer this question, we use Examples \ref{exam-p0-1} and \ref{exam-p0-2}. It has been proved that if $\mathscr{A}$ is given by Example \ref{exam-p0-1} or Example \ref{exam-p0-2}, then $\mathscr{A}\in Q$ but $\mathscr{A}\not\in R_0$. In addition,
\begin{itemize}
\item if $\mathscr{A}$ is given by Example \ref{exam-p0-1}, then $x^{T}\mathscr{A}x^{3}=x_2^4\geq0$ for any $x\in \mathbb{R}^2_+$ with $x\neq0$; and
\item if $\mathscr{A}$ is given by Example \ref{exam-p0-2}, then $x^{T}\mathscr{A}x^{2}=x_2^3\geq0$ for any $x\in \mathbb{R}^2_+$ with $x\neq0$,
\end{itemize}
which imply that if $\mathscr{A}$ is given by Example \ref{exam-p0-1} or Example \ref{exam-p0-2}, then $\mathscr{A}$ is copositive. Thus, the result of Proposition \ref{pro-p-22} cannot be extended to the case of tensor. Moreover, by using Example \ref{exam-q-l1-2} we can also obtain the above result.

\section{Conclusions}

In this paper, we studied $Q$-tensor and gave a sufficient condition to judge whether a tensor is a $Q$-tensor or not. In order to extend Agangic-Cottle's result to the tensor space, we introduced the concept of $SP_0$-tensor and showed that within the class of $SP_0$-tensors, four classes of tensors, i.e., $R_0$-tensors, $R$-tensors, $ER$-tensors and $Q$-tensors, are all equivalent. We clarified that the above equivalence does not hold if $SP_0$-tensors is replaced by $P_0$-tensors or $P_0^\prime$-tensors by constructing two examples; and discussed the relationships among three classes of $P_0$-type tensors, i.e., $P_0$-tensors, $P_0^\prime$-tensors and $SP_0$-tensors. In order to extend Danao's result to the tensor space, we showed that the above equivalence holds within the class of nonnegative tensors. We also discussed the relationship between the class of $SP_0$-tensors and the class of nonnegative tensors. Moreover, by using several examples we also showed that three famous results, related to $Q$-matrices, semi-positive tensors and copositive tensors, cannot be extended to the tensor space. Since $Q$-matrix plays an important role in the field of structured matrices and the theory of complementarity problems, fruitful results related to $Q$-matrices have been obtained in the literature. We believe more results related to $Q$-matrices can be clarified whether they can be extended to the case of tensor or not.

%
%

\end{document}